\documentclass[letterpaper, 10 pt, conference]{ieeeconf}  

\IEEEoverridecommandlockouts                              

\usepackage{graphics} 
\usepackage{epsfig} 
\usepackage{times} 
\usepackage{amsmath} 
\usepackage{amssymb}  
\usepackage{bbm}
\usepackage{algorithm}
\usepackage{cite}
\usepackage[noend]{algpseudocode}

\usepackage{dsfont}

\newtheorem{assumption}{Assumption}

\newtheorem{lemma}{Lemma}

\allowdisplaybreaks[1]  

\title{\LARGE \bf
Hierarchical Optimal Power Flow with Improved Gradient Evaluation}

\author{Heng Liang,$^{1}$ Xinyang Zhou,$^{2}$ and Changhong Zhao$^{1}$
\thanks{This work was supported by Hong Kong Research Grants Council through ECS Award No. 24210220 and by CUHK through faculty startup fund.}
\thanks{$^{1}$H. Liang and C. Zhao are with the Department of Information Engineering, the Chinese University of Hong Kong, New Territories, HKSAR, China. 
        {\tt\small \{lh021, chzhao\}@ie.cuhk.edu.hk}}%
\thanks{$^{2}$X. Zhou is with the National Renewable Energy Laboratory, Golden, CO, USA.
        {\tt\small xinyang.zhou@nrel.gov}}%
\thanks{\copyright 2022 IEEE. Personal use of this material is permitted. Permission from IEEE must be obtained for all other uses, in any current or future media, including reprinting/republishing this material for advertising or promotional purposes, creating new collective works, for resale or redistribution to servers or lists, or reuse of any copyrighted component of this work in other works. This paper is accepted for presentation at the American Control Conference, Atlanta, GA, USA, 2022.}
}

\begin{document}

\maketitle
\thispagestyle{empty}
\pagestyle{empty}

\begin{abstract}

Existing algorithms to solve alternating-current optimal power flow (AC-OPF) often exploit linear approximations to simplify system models and accelerate computations. In this paper, we improve a recent hierarchical OPF algorithm, which rested on primal-dual gradients evaluated in a linearized distribution power flow model. Specifically, we identify a risk of voltage violation arising from the model linearization, and propose a more accurate gradient evaluation method to eliminate that risk. 
We further develop a hierarchical primal-dual algorithm to solve OPF based on the proposed gradient evaluation method. Numerical results on IEEE networks show that our algorithm can enhance voltage safety with satisfactory computational efficiency.  
\end{abstract}

\section{Introduction}

Optimal power flow (OPF) is a fundamental optimization problem that aims to find a cost-minimizing operating point subject to the physical laws and safety limits of a power network. 
The power sector is witnessing a rising need for fast and scalable OPF solvers, which can instruct a large network of controllable units (smart appliances, electric vehicles, energy storage devices, etc.) to make timely response to the growing variations of wind and solar generations.
Such a need is especially urgent in power distribution networks, where massive renewable energy sources and controllable units are being deployed.
However, distribution networks are also where the speed and scalability requirements are most challenging, as the high resistance-to-reactance ratios of distribution lines necessitate nonlinear and nonconvex alternating-current (AC) OPF rather than its simple direct-current approximation.

Numerous efforts have been made to overcome the computational challenges to AC-OPF. Many of them conducted convex relaxations, convex inner approximations, or linearizations of AC power flow; comprehensive reviews were provided in \cite{low2014convex, molzahn2019survey}.
Meanwhile, distributed OPF algorithms were developed and shown to be more scalable in terms of computation and communication and more robust to single-point failure, compared to their centralized counterparts \cite{dall2013distributed, erseghe2014distributed, zhang2014optimal, peng2016distributed}.
To further reduce computational efforts associated with solving AC power flow, some OPF algorithms were implemented by iteratively actuating the power system with intermediate decisions and updating the decisions based on system feedback \cite{bolognani2014distributed, gan2016online, bernstein2019real}. 

From the vast literature, we bring to attention a hierarchical distributed primal-dual gradient algorithm \cite{zhou2019hierarchical} (and its extension to multi-phase networks \cite{zhou2019accelerated}). 
It leveraged the radial structure of distribution networks to avoid repetitive computation and communication, and thus significantly accelerated large-scale OPF computations.
However, the algorithm in \cite{zhou2019hierarchical} 
derived approximate gradients from the linearized distribution power flow model in \cite{baran1989optimalC, farivar2013branch}. Such linearization may cause the solver to optimistically estimate nodal voltages to be safe, while they actually already exceed safety limits. The consequent risk of voltage violation will be revealed later in this paper. 

To prevent such violation, we develop an improved gradient evaluation method by approximating the partial derivatives of the quadratic terms associated with line currents and power losses.
Our analysis shows that with moderate extra computations, the proposed method returns more accurate gradient estimations than \cite{zhou2019hierarchical}.  
Meanwhile, the proposed method preserves the structure in \cite{zhou2019hierarchical} for gradient evaluation, and thus enables us to develop an improved hierarchical OPF algorithm. 
Numerical experiments on IEEE 37-node and 123-node networks demonstrate enhanced safety of voltage regulation achieved by the proposed algorithm with a moderate increase in computation time, compared to the previous method based on the linearized model.  

The rest of this paper is organized as follows. Section \ref{sec:model} introduces the distribution network model, the OPF problem, and a primal-dual algorithm to solve it. Section \ref{sec:method} motivates and elaborates the improved gradient evaluation method. Section \ref{sec:algorithm} presents the hierarchical OPF algorithm based on the improved gradient evaluation. 
Section \ref{sec:simulation} reports the numerical experiments. Section \ref{sec:conclusion} concludes this paper.

\section{Modeling and Preliminary Algorithm}\label{sec:model}

\subsection{Power network model and OPF problem}

We model a single-phase equivalent power distribution network as a directed \emph{tree} graph $\mathcal{T}:=\{\mathcal{N}^{+}, \mathcal{E}\}$, where $\mathcal{N}^{+}=\mathcal{N} \cup\{0\}$, with $0$ indexing the root node (the substation, also known as the slack bus), and $\mathcal{N}=\{1,...,N\}$ indexing other nodes. 
The set $\mathcal{E}$ collects the ordered pairs of nodes representing the lines in the network. 
We define the directions of lines as pointing away from the root, e.g., in line $(i,j)$, node $i$ is closer to the root than node $j$.
Let $p_{i}$ and $q_{i}$ denote the net active and reactive power injections (i.e., power supply minus consumption) and $v_{i}$ denote the squared voltage magnitude at each node $i\in \mathcal{N}^+$. 
Let $\ell_{ij}$ denote the squared current magnitude, $z_{ij}=r_{ij}+\boldsymbol{i}x_{ij}$ denote the series impedance, and $P_{ij}$ and $Q_{ij}$ denote the sending-end active and reactive power on each line $(i,j)\in\mathcal{E}$.

We adopt the distribution power flow model \cite{baran1989optimalC, farivar2013branch}:
\begin{subequations}\label{powerequation}
\begin{alignat}{2}
P_{i j} &=-p_{j}+\sum_{k:(j, k) \in \mathcal{E}} P_{j k}+r_{i j} \ell_{i j},~\forall j\in\mathcal{N} \label{powerequation:1}\\
Q_{i j} &=-q_{j}+\sum_{k:(j, k) \in \mathcal{E}} Q_{j k}+x_{i j} \ell_{i j}, ~\forall j\in\mathcal{N} \label{powerequation:2}\\
v_{j} &=v_{i}\!-\!2\left(r_{i j} P_{i j}\!+\!x_{i j} Q_{i j}\right)\!+\!\left|z_{i j}\right|^{2} \ell_{i j}, \forall (i,j)\in\mathcal{E}  \label{powerequation:3}\\
\ell_{i j} v_{i} &=P_{i j}^{2}+Q_{i j}^{2},~\forall (i,j)\in\mathcal{E}. \label{powerequation:4}
\end{alignat}
\end{subequations}
Assume the power injections $(p_i,q_i)$ to be controllable within a given operating region:
\begin{alignat}{2} \label{region}
\mathcal{Y}_{i}=\left\{\left(p_{i}, q_{i}\right) \mid \underline{p}_{i} \leqslant p_{i} \leqslant \bar{p}_{i}, \underline{q}_{i} \leqslant q_{i} \leqslant \bar{q}_{i}\right\},~\forall i\in\mathcal{N}.
\end{alignat}
Define vectors $\boldsymbol{p}:=\left[p_{1}, \ldots, p_{N}\right]^{\top}$, $\boldsymbol{q}:=\left[q_{1}, \ldots, q_{N}\right]^{\top}$, $\boldsymbol{v}:=\left[v_{1}, \ldots, v_{N}\right]^{\top} \in \mathbb{R}^{N}$.
Following \cite{gan2016online}, we express the voltage vector as a function of power injections, i.e., $\boldsymbol{v}(\boldsymbol{p},\boldsymbol{q})$, implicitly defined by \eqref{powerequation}. 
Consider the following OPF problem:
\begin{subequations} \label{OPF}
\begin{alignat}{2}
\min _{\boldsymbol{p}, \boldsymbol{q}} & \sum_{i \in \mathcal{N}} f_{i}\left(p_{i}, q_{i}\right) \label{OPF:obj}\\
\text { s.t. } & \underline{\boldsymbol{v}} \leqslant \boldsymbol{v}(\boldsymbol{p}, \boldsymbol{q}) \leqslant \overline{\boldsymbol{v}} \label{OPF:v}\\
& \left(p_{i}, q_{i}\right) \in \mathcal{Y}_{i}, ~\forall i \in \mathcal{N}
\end{alignat}
\end{subequations}
where $f_{i}$ is a \emph{strongly convex} cost function for the control of power injections at node $i\in\mathcal{N}$, inequality \eqref{OPF:v} is element-wise, given constant squared voltage limits $\underline{\boldsymbol{v}}$ and $\overline{\boldsymbol{v}}$.

\subsection{Primal-dual gradient algorithm}

Let $\underline{\boldsymbol{\mu}}$ and $\overline{\boldsymbol{\mu}}$ be the dual variables associated with \eqref{OPF:v} to penalize the violation of voltage lower and upper bounds, respectively. 
To design a convergent primal-dual algorithm, we consider the regularized Lagrangian of \eqref{OPF}:
\begin{alignat}{2}\label{RLagarangian} \nonumber
\mathcal{L}_{\epsilon}(\boldsymbol{p}, \boldsymbol{q} ; \overline{\boldsymbol{\mu}}, \underline{\boldsymbol{\mu}})= \sum_{i \in \mathcal{N}} f_{i}\left(p_{i}, q_{i}\right)\qquad\qquad\\ 
+\underline{\boldsymbol{\mu}}^{\top}(\underline{\boldsymbol{v}}-\boldsymbol{v}(\boldsymbol{p}, \boldsymbol{q}))+\overline{\boldsymbol { \mu }}^{\top}(\boldsymbol{v}(\boldsymbol{p}, \boldsymbol{q})-\overline{\boldsymbol{v}})-\frac{\epsilon}{2}\|\boldsymbol{\mu}\|_{2}^{2}
\end{alignat}
where $\boldsymbol{\mu}=\left(\underline{\boldsymbol{\mu}},\overline{\boldsymbol{\mu}}\right)$ and $\epsilon>0$ is a constant factor for the regularization term. 
Function \eqref{RLagarangian} is strongly convex in primal variables $(\boldsymbol{p}, \boldsymbol{q})$ and strongly concave in dual variables $\boldsymbol{\mu}$. The saddle point of \eqref{RLagarangian} serves as an approximate solution to \eqref{OPF}, with a bounded error that is linear in $\epsilon$ \cite{zhou2019hierarchical,koshal2011multiuser}.
The following primal-dual gradient algorithm has been commonly used to iteratively approach such a saddle point:
\begin{subequations} \label{iter}
\begin{alignat}{2} \nonumber
\boldsymbol{u}(t+1)=\bigg[\boldsymbol{u}(t)-\sigma_u\ \left(\frac{\partial f\left(\boldsymbol{u}(t)\right)}{\partial \boldsymbol{u}}\right)^{\top}   \nonumber\\
-\sigma_u  \left(\frac{\partial \boldsymbol{v}\left(\boldsymbol{u}(t)\right)}{\partial \boldsymbol{u} }\right)^{\top} \left(\overline{\boldsymbol{\mu}}(t)-\underline{\boldsymbol{\mu}}(t)\right)\bigg]_{\mathcal{Y}} \label{iter:u}\\ 
\underline{\boldsymbol{\mu}}(t+1)=\left[\underline{\boldsymbol{\mu}}(t)+\sigma_{\mu}\left(\underline{\boldsymbol{v}}-\boldsymbol{v}(t)-\epsilon \underline{\boldsymbol{\mu}}(t)\right)\right]_{+}\label{iter:underline}\\
\overline{\boldsymbol{\mu}}(t+1)=\left[\overline{\boldsymbol{\mu}}(t)+\sigma_{\mu}\left(\boldsymbol{v}(t)-\overline{\boldsymbol{v}}-\epsilon \overline{\boldsymbol{\mu}}(t)\right)\right]_{+}\label{iter:overline}
\end{alignat}
\end{subequations}
where $\boldsymbol{u}:=(\boldsymbol{p},\boldsymbol{q})$ collects the controllable power injections; $f(\boldsymbol{u}):=\sum_{i\in\mathcal{N}} f_i(p_i,q_i)$ is the objective \eqref{OPF:obj}; $\boldsymbol{v}(t)=\boldsymbol{v}\left(\boldsymbol{u}(t)\right)$ both denote the voltage profile under power injections $\boldsymbol{u}(t)$ at time step $t$. Note that $\boldsymbol{v}\left(\boldsymbol{u}(t)\right)$ does not have to be calculated by \eqref{powerequation} or any other mathematical model of power flow; instead, it can be measured from the power network after actuating the network with the intermediate decisions $\boldsymbol{u}(t)$ during the process \eqref{iter}. This \emph{feedback-based} implementation has been widely adopted to compensate for likely modeling errors \cite{bolognani2014distributed, gan2016online, bernstein2019real, zhou2019accelerated}. 
The subscripts $\mathcal{Y}$ and $+$ represent the projections onto $\prod_{i\in\mathcal{N}} \mathcal{Y}_{i}$ and the nonnegative orthant, respectively. Without loss of generality, we use constant step sizes $\sigma_u, \sigma_{\mu}>0$ to update the primal and dual variables, respectively.        

The gradient $\partial{f}/\partial{\boldsymbol{u}}$ in \eqref{iter:u} can be easily computed, and the voltage $\boldsymbol{v}(t)$ in \eqref{iter:underline}--\eqref{iter:overline} can be conveniently measured, both locally at each node $i\in\mathcal{N}$. 
Therefore, the remaining key challenge is to compute the gradient $\partial{\boldsymbol{v}}/\partial{\boldsymbol{u}}$ that potentially couples all the nodes in the power network. This will be addressed in the next section.

\section{Improved Gradient Evaluation}
\label{sec:method}

\subsection{Prior work and motivation for improvement}

In power flow equations \eqref{powerequation}, the nonlinear and implicit dependence of voltage $\boldsymbol{v}$ on power injections $\boldsymbol{u}$ makes it difficult to quickly and precisely compute the gradient $\partial{\boldsymbol{v}}/\partial{\boldsymbol{u}}$. 
Reference \cite{gan2016online} proposed a backward-forward sweep method for gradient calculation, which is computationally inefficient in large networks. In contrast, prior work \cite{zhou2019hierarchical} considered the linearized power flow model in \cite{baran1989optimalC, farivar2013branch}: 
\begin{subequations}\label{simppowerequation}
\begin{alignat}{2}
\hat{P}_{i j} &=-p_{j}+\sum_{k:(j, k) \in \mathcal{E}} \hat{P}_{j k},~\forall j\in\mathcal{N}  \label{simppowerequation:1}\\
\hat{Q}_{i j} &=-q_{j}+\sum_{k:(j, k) \in \mathcal{E}} \hat{Q}_{j k},~\forall j\in\mathcal{N}  \label{simppowerequation:2}\\
\hat{v}_{j} &=\hat{v}_{i}-2\left(r_{i j} \hat{P}_{i j}+x_{i j} \hat{Q}_{i j}\right),~\forall (i,j)\in\mathcal{E} \label{simppowerequation:3}
\end{alignat}
\end{subequations}
which yields the following approximation of $\partial{\boldsymbol{v}}/\partial{\boldsymbol{u}}$:
\begin{subequations}\label{linear-gradient}
\begin{alignat}{2}
\frac{\partial \hat{v}_{j}(\boldsymbol{p},\boldsymbol{q})}{\partial{p_{h}}} = R_{j h}:=\sum_{(i, k) \in \mathbb{P}_{j \wedge h}} 2\cdot r_{i k},~\forall j,h \in \mathcal{N}  \label{linear-gradient:p}\\
\frac{\partial \hat{v}_{j}(\boldsymbol{p},\boldsymbol{q})}{\partial{q_{h}}} = X_{j h}:=\sum_{(i, k) \in \mathbb{P}_{j \wedge h}} 2\cdot x_{i k},~\forall j,h \in \mathcal{N} \label{linear-gradient:q}
\end{alignat}
\end{subequations}
where $\mathbb{P}_{j \wedge h}$ denotes the common part of the unique paths from nodes $j$ and $h$ back to the root node.


We illustrate the error of linear model \eqref{simppowerequation} by introducing the following \emph{Lemma}.
\begin{lemma}\label{lemma::1}
Let $g_{\xi}$ and $\eta_{\xi}$ denote the active and reactive power loss downstream of node $\xi\in \mathcal{N}$. The voltage error induced by \eqref{simppowerequation} for all $h\in \mathcal{N}$ is:
\begin{alignat}{2} \label{verror}
\hat{v}_{h}-v_{h}= \sum_{(\zeta,\xi)\in \mathbb{P}_{h}}\left(2\left(r_{\zeta\xi}g_{\xi}+x_{\zeta\xi}\eta_{\xi}\right)+(r_{\zeta\xi}^2+x_{\zeta\xi}^2)\ell_{\zeta\xi}\right),
\end{alignat}
where the expressions of $g_{\xi}$ and $\eta_{\xi}$ are:
\begin{alignat}{2} \label{geta}
g_{\xi}=\sum_{(\alpha,\beta)\in \operatorname{down}(\xi)}r_{\alpha\beta}\ell_{\alpha\beta},\quad \eta_{\xi}=\sum_{(\alpha,\beta)\in \operatorname{down}(\xi)}x_{\alpha\beta}\ell_{\alpha\beta}. \nonumber
\end{alignat}
\end{lemma}

\proof See Appendix.


\emph{Lemma \ref{lemma::1}} reveals that the voltage computed with the linear model \eqref{simppowerequation} is higher than that under the accurate nonlinear model \eqref{powerequation}. Such errors would accumulate as we trace the nodes further away from the root. 
Therefore, under the power injections $(\boldsymbol{p}, \boldsymbol{q})$ solved by \eqref{iter} with the linear model \eqref{simppowerequation} and the approximate gradient \eqref{linear-gradient}, even though the model optimistically estimates the voltages to be safe, the actual voltages may already drop below their lower bounds. The need to prevent such voltage violation motivates us to develop an improved gradient evaluation method, as elaborated below.

\subsection{Improved gradient evaluation}

Consider an arbitrary node $h\in\mathcal{N}$, and $u_h := (p_h,q_h)$. We take the partial derivatives of the variables in the nonlinear power flow model \eqref{powerequation}:
\begin{subequations}\label{derive}
\begin{alignat}{2}
\frac{\partial P_{i j}}{\partial u_{h}} =&-\frac{\partial p_{j}}{\partial u_{h}}+\sum_{k:(j, k) \in \mathcal{E}} \frac{\partial P_{j k}}{\partial u_{h}}+r_{i j} \frac{\partial \ell_{i j}}{\partial u_{h}}, ~\forall j\in\mathcal{N} \\
\frac{\partial Q_{i j}}{\partial u_{h}} =&-\frac{\partial q_{j}}{\partial u_{h}}+\sum_{k:(j, k) \in \mathcal{E}} \frac{\partial Q_{j k}}{\partial u_{h}}+x_{i j} \frac{\partial \ell_{i j}}{\partial u_{h}}, \forall j\in\mathcal{N} \\
\frac{\partial v_{j}}{\partial u_{h}} =&\frac{\partial v_{i}}{\partial u_{h}}\!-\!2\left(r_{i j}\frac{\partial P_{i j}}{\partial u_{h}} \!+\!x_{i j}\frac{\partial Q_{i j}}{\partial u_{h}} \right)  \!+\!\left|z_{i j}\right|^{2} \frac{\partial \ell_{i j}}{\partial u_{h}} \label{derive:c}\\
\frac{\partial\ell_{i j}}{\partial u_{h}} =&\frac{2P_{i j}}{v_{i}}\frac{\partial P_{ij}}{\partial u_{h}}+\frac{2Q_{i j}}{v_{i}}\frac{\partial Q_{ij}}{\partial u_{h}}-\frac{\ell_{ij}}{v_{i}}\frac{\partial v_{i}}{\partial u_{h}},~\forall (i,j)\in\mathcal{E}. \label{derive:d}
\end{alignat}
\end{subequations}

The exact partial derivatives are hard to solve from \eqref{derive} due to their complex interdependence. To facilitate computation, we first simplify \eqref{derive:d} by replacing the partial derivatives with their approximations under the linear model \eqref{simppowerequation}: 
\begin{alignat}{2} \label{deriveell}
\frac{\partial\hat{\ell}_{i j}}{\partial u_{h}} &=\frac{2P_{i j}}{v_{i}}\frac{\partial \hat{P}_{ij}}{\partial u_{h}}+\frac{2Q_{i j}}{v_{i}}\frac{\partial \hat{Q}_{ij}}{\partial u_{h}}-\frac{\ell_{ij}}{v_{i}}\frac{\partial \hat{v}_{i}}{\partial u_{h}}.
\end{alignat}
We recall the following results from \cite{gan2016online} for all the nodes $h\in\mathcal{N}$ and lines $(i,j)\in\mathcal{E}$:
\begin{subequations} \label{lineargradient-PQ}
\begin{alignat}{2}
\frac{\partial{\hat{P}_{i j}}}{\partial p_h} &=-\mathds{1}\left(j \in \mathbb{P}_{h}\right),\quad & \frac{\partial{\hat{P}_{i j}}}{\partial q_h} &=0 \\
\frac{\partial{\hat{Q}_{i j}}}{\partial q_h} &=-\mathds{1}\left(j \in \mathbb{P}_{h}\right), \quad & \frac{\partial{\hat{Q}_{i j}}}{\partial p_h} &=0 
\end{alignat}
\end{subequations}
where $\mathds{1}(j \in \mathbb{P}_{h})$ is an indicator that equals 1 if node $j$ lies on the unique path from node $h$ to the root, and 0 otherwise. 
By \eqref{lineargradient-PQ} and \eqref{linear-gradient}, we can calculate $\partial\hat{\boldsymbol{\ell}}/\partial {\boldsymbol{u}}$ in \eqref{deriveell} as follows:
\begin{subequations}\label{deriveellresult}
\begin{alignat}{2}
\frac{\partial\hat{\ell}_{i j}}{\partial p_{h}} &=-\frac{1}{v_{i}}\left(2 P_{i j}\cdot \mathds{1}(j \in \mathbb{P}_{h})+\ell_{ij}R_{i  h}\right)\\
 \frac{\partial\hat{\ell}_{i j}}{\partial q_{h}} &=-\frac{1}{v_{i}}\left(2 Q_{i j}\cdot \mathds{1}(j \in \mathbb{P}_{h})+\ell_{ij}X_{i  h}\right). 
\end{alignat}
\end{subequations}

The results in \eqref{linear-gradient}, \eqref{lineargradient-PQ}, \eqref{deriveellresult} can replace the partial derivatives on the right-hand side (RHS) of \eqref{derive:c} to obtain an \emph{improved gradient evaluation} $\partial{\boldsymbol{v}}/\partial{\boldsymbol{u}}$ as: 
\begin{subequations}\label{derivesecond}
\begin{alignat}{2}
&\frac{\partial v_{j}}{\partial p_{h}} =\frac{\partial \hat{v}_{i}}{\partial p_{h}}-2\left(r_{i j}\frac{\partial \hat{P}_{i j}}{\partial p_{h}} +x_{i j}\frac{\partial \hat{Q}_{i j}}{\partial p_{h}} \right)+\left|z_{i j}\right|^{2} \frac{\partial \hat{\ell}_{i j}}{\partial p_{h}} \nonumber \\
& =\left(1\!-\!\frac{\left|z_{ij}\right|^{2}\ell_{ij}}{v_{i}}\right)R_{ih}+2\left(r_{ij}\!-\!\frac{\left|z_{ij}\right|^{2}P_{ij}}{v_{i}}\right)\cdot \mathds{1}(j \in \mathbb{P}_{h})\label{derivesecond:p} \\
&\frac{\partial v_{j}}{\partial q_{h}} =\frac{\partial \hat{v}_{i}}{\partial q_{h}}-2\left(r_{i j}\frac{\partial \hat{P}_{i j}}{\partial q_{h}} +x_{i j}\frac{\partial \hat{Q}_{i j}}{\partial q_{h}} \right)+\left|z_{i j}\right|^{2} \frac{\partial \hat{\ell}_{i j}}{\partial q_{h}} \nonumber \\
& =\left(1\!-\!\frac{\left|z_{ij}\right|^{2}\ell_{ij}}{v_{i}}\right)X_{ih}+2\left(x_{ij}\!-\!\frac{\left|z_{ij}\right|^{2}Q_{ij}}{v_{i}}\right)\cdot \mathds{1}(j \in \mathbb{P}_{h}) \label{derivesecond:q}
\end{alignat}
\end{subequations}
where node $i$ is the unique parent of node $j$ in the directed tree network.  
By \eqref{linear-gradient:p}, there is $R_{jh} = R_{ih}+ 2r_{ij}\mathds{1}(j \in \mathbb{P}_{h})$, which converts the result in \eqref{derivesecond:p} into:
\begin{eqnarray}\label{derivesecond:rewrite:p}
\frac{\partial v_{j}}{\partial p_{h}} = R_{jh}  -\frac{\left|z_{ij}\right|^{2}\ell_{ij}}{v_{i}} R_{ih}-\frac{2\left|z_{ij}\right|^{2}P_{ij}}{v_{i}} \mathds{1}(j \in \mathbb{P}_{h}). 
\end{eqnarray}
The first term on the RHS of \eqref{derivesecond:rewrite:p} is the same as \eqref{linear-gradient:p} derived from the linear model \eqref{simppowerequation}, while the second and third terms compensate for the effects of the quadratic terms that were neglected from \eqref{simppowerequation}. 
The improved partial derivatives \eqref{derivesecond:q} with respect to reactive power injections $q_h$ have the same structure as \eqref{derivesecond:p}, and thus we will not repeat.

\section{Hierarchical OPF Algorithm}\label{sec:algorithm}

Based on the primal-dual framework \eqref{iter} and the improved gradient evaluation \eqref{derivesecond}, we design a scalable hierarchical algorithm to solve the OPF problem \eqref{OPF}. The proposed algorithm utilizes the subtree structure of radial distribution networks discovered in \cite{zhou2019accelerated, zhou2019hierarchical}, while achieving more accurate and safer voltage regulation. 

The distribution tree network $\mathcal{T}:=\{\mathcal{N}^{+}, \mathcal{E}\}$ is composed of subtrees  
$\mathcal{T}_{k}=\left\{\mathcal{N}_{k}, \mathcal{E}_{k}\right\}$ indexed by $k \in \mathcal{K}=\{1, \ldots, K\}$ and a set of nodes $\mathcal{N}_{0}$ that are not clustered into any subtree. Let $n_{k}^{0}$ denote the root node of subtree $k$, which is the node in $\mathcal{N}_{k}$ that is nearest to the root node of the whole network.
Given a distribution network, the division of subtrees and unclustered nodes may not be unique, but we assume it always satisfies the following conditions.
\begin{assumption}\label{ass:nonoverlapping}
The subtrees are \emph{non-overlapping}, i.e., $\mathcal{N}_{k_1} \cap \mathcal{N}_{k_2}=\emptyset$ for any $k_1, k_2\in\mathcal{K}$, $k_1\neq k_2$. 
\end{assumption}
\begin{assumption}\label{ass:notonpath}
For any subtree root node $n_{k}^{0}$, $k\in \mathcal{K}$ or any unclustered node $n\in \mathcal{N}_0$, its path to the network root node only goes through a subset of nodes in $\mathcal{N}_0$, but not any node in another subtree. 
\end{assumption}

\begin{figure} 
    \centering
    \includegraphics[width = 0.49\columnwidth]{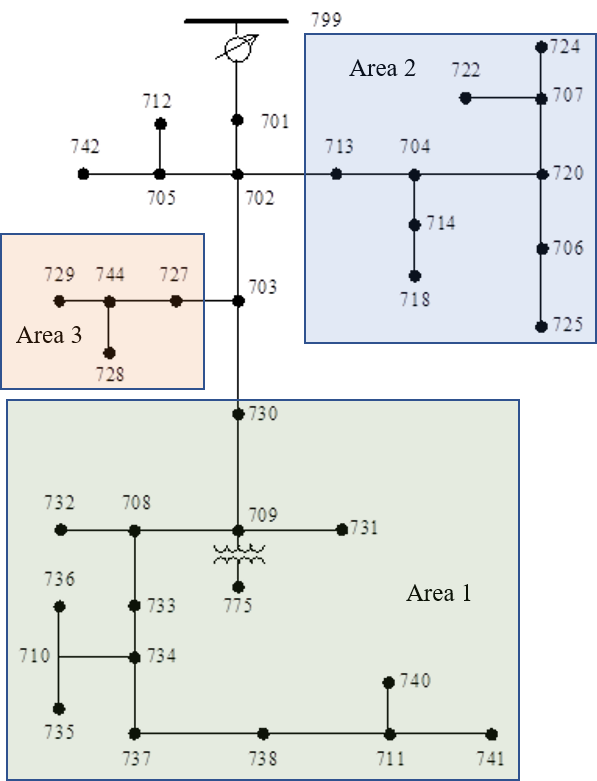}
    \hfil
    \includegraphics[width = 0.49\columnwidth]{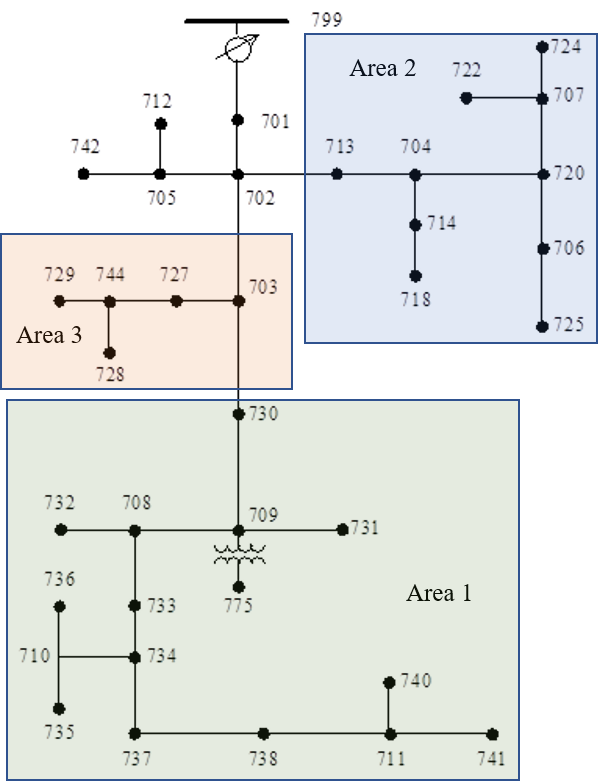}
    \caption{The clustering of the IEEE 37-node network on the left satisfies Assumption \ref{ass:notonpath}. The one on the right does not, because the path from a subtree root node 730 to the network root 799 goes through another subtree.}
    \label{ieee37}
\end{figure}

The left subfigure of Figure \ref{ieee37} shows a clustering of the IEEE 37-node network that satisfies Assumption \ref{ass:notonpath}, while the one in the right subfigure does not.


The distribution system operator plays a role of central controller (CC). It is sufficient for the CC to know the structure and parameters of the \emph{backbone network}, which only connects the unclustered nodes $\mathcal{N}_{0}$ and the subtree root nodes $\left\{n_{k}^0,~ k\in \mathcal{K}\right\}$. The CC also maintains two-way communications to receive/send/relay information from/to/between the backbone network nodes.
Each subtree network $k\in\mathcal{K}$ is known to and managed by the $k$-th regional controller (RC $k$), which communicates with the nodes in subtree $k$, the parent node of the subtree root $n_k^0$, and the CC.

With the settings above, we now explain the hierarchical structure underlying the vector $\left(\frac{\partial \boldsymbol{v}}{\partial \boldsymbol{u} }\right)^{\top} \left(\overline{\boldsymbol{\mu}}-\underline{\boldsymbol{\mu}}\right)$ in \eqref{iter:u}.
Without loss of generality, we only consider the element of this vector corresponding to $p_h$ for a particular node $h\in\mathcal{N}$, which can be discussed in two cases below. The element corresponding to $q_h$ can be calculated in the same manner. 

\textit{Case 1:} If $h\in \mathcal{N}_k$ is in a subtree $k\in\mathcal{K}$, we have:
\begin{eqnarray}
 &&\quad \sum_{j \in \mathcal{N}}\frac{\partial v_{j}}{\partial p_{h}}\cdot\left(\overline{\mu}_{j}-\underline{\mu}_{j}\right) =\sum_{j \in \mathcal{N}_{k}}\frac{\partial v_{j}}{\partial p_{h}}\cdot\left(\overline{\mu}_{j}-\underline{\mu}_{j}\right) \nonumber\\
 &&\qquad +\sum_{k' \in \mathcal{K}\backslash\{k\}}~\sum_{j \in \mathcal{N}_{k'}}\frac{\partial v_{j}}{\partial p_{h}}\cdot\left(\overline{\mu}_{j}-\underline{\mu}_{j}\right)
 \nonumber 
 \\ && \qquad +\sum_{j \in \mathcal{N}_{0}}\frac{\partial v_{j}}{\partial p_{h}}\cdot\left(\overline{\mu}_{j}-\underline{\mu}_{j}\right).\label{eq:decomposed_dvdp}
\end{eqnarray}

The first term on the RHS of \eqref{eq:decomposed_dvdp} sums over nodes $j$ in the same subtree $k$ as node $h$. It can be calculated using \eqref{derivesecond:p} by RC $k$, and then sent to node $h$ which requires this information to compute \eqref{eq:decomposed_dvdp} and carry out \eqref{iter:u}. In particular, the calculation of $\partial v_j/\partial p_h$ requires RC $k$ to receive parameters $R_{ih}$, $z_{ij}$ and measurements $\ell_{ij}$, $P_{ij}$, $v_i$ from line $(i,j)$ that connects node $j$ to its parent node $i$.  


The second term on the RHS of \eqref{eq:decomposed_dvdp} sums over nodes $j$ in all the subtrees $k'$ except subtree $k$ that hosts node $h$. By Assumption \ref{ass:notonpath}, there must be $j\notin \mathbb{P}_{h}$, which simplifies \eqref{derivesecond:p} and leads to:
\begin{eqnarray} 
\sum_{k' \in \mathcal{K}\backslash\{k\}}~\sum_{j \in \mathcal{N}_{k'}}\frac{\partial v_{j}}{\partial p_{h}}\cdot\left(\overline{\mu}_{j}-\underline{\mu}_{j}\right)\qquad\qquad\qquad&& \nonumber \\
= \sum_{k' \in \mathcal{K}\backslash\{k\}}
R_{n_{k}^{0}n_{k'}^{0}}\sum_{j \in \mathcal{N}_{k'}}\left(1\!-\!\frac{\left|z_{ij}\right|^{2}\ell_{ij}}{v_{i}}\right)\cdot\left(\overline{\mu}_{j}\!-\!\underline{\mu}_{j}\right).&& \label{NkNt}
\end{eqnarray}
In \eqref{NkNt}, the sum over $j\in \mathcal{N}_{k'}$ can be calculated by RC $k'$. The CC receives such sums from all the RCs $k'\in \mathcal{K}\backslash\{k\}$, adds them up after weighting by $R_{n_{k}^{0}n_{k'}^{0}}$, and sends the result to RC $k$. Upon receiving this result, RC $k$ broadcasts it to all the nodes in subtree $k$, including node $h$ which requires this information to compute \eqref{eq:decomposed_dvdp} and carry out \eqref{iter:u}. 

The third term on the RHS of \eqref{eq:decomposed_dvdp} sums over all the unclustered nodes $j\in\mathcal{N}_0$. It can be calculated using \eqref{derivesecond:p} by the CC and then relayed by RC $k$ to node $h$. In particular, $R_{ih} = R_{i n_k^0}$ and $\mathds{1}(j\in\mathbb{P}_{h}) = \mathds{1}(j\in\mathbb{P}_{n_k^0})$ are known to the CC for the calculation of \eqref{derivesecond:p}.


\textit{Case 2:} If $h\in \mathcal{N}_0$ is an unclustered node, we have: 
\begin{eqnarray}
 \sum_{j \in \mathcal{N}}\frac{\partial v_{j}}{\partial p_{h}}\cdot\left(\overline{\mu}_{j}-\underline{\mu}_{j}\right) &=& \sum_{k \in \mathcal{K}}~\sum_{j \in \mathcal{N}_{k}}\frac{\partial v_{j}}{\partial p_{h}}\cdot\left(\overline{\mu}_{j}-\underline{\mu}_{j}\right)
 \nonumber 
 \\ && +\sum_{j \in \mathcal{N}_{0}}\frac{\partial v_{j}}{\partial p_{h}}\cdot\left(\overline{\mu}_{j}-\underline{\mu}_{j}\right). \label{eq:decomposed_dvdp_unclustered}
\end{eqnarray}

The first term on the RHS of \eqref{eq:decomposed_dvdp_unclustered} sums over nodes $j$ in all the subtrees $k\in\mathcal{K}$. By Assumption \ref{ass:notonpath}, there must be $j\notin \mathbb{P}_{h}$, which simplifies \eqref{derivesecond:p} and leads to:
\begin{eqnarray} 
\sum_{k \in \mathcal{K}}~\sum_{j \in \mathcal{N}_{k}}\frac{\partial v_{j}}{\partial p_{h}}\cdot\left(\overline{\mu}_{j}-\underline{\mu}_{j}\right)\qquad\qquad\qquad&& \nonumber \\
= \sum_{k \in \mathcal{K}}
R_{h n_{k}^{0}}\sum_{j \in \mathcal{N}_{k}}\left(1-\frac{\left|z_{ij}\right|^{2}\ell_{ij}}{v_{i}}\right)\cdot\left(\overline{\mu}_{j}-\underline{\mu}_{j}\right).&& \label{NkNt_unclustered}
\end{eqnarray}
In \eqref{NkNt_unclustered}, the sum over $j\in \mathcal{N}_{k}$ can be calculated by RC $k$. The CC receives such sums from all the RCs $k\in \mathcal{K}$, adds them up after weighting by $R_{h n_{k}^{0}}$, and sends the result to node $h$ for subsequent computations of \eqref{eq:decomposed_dvdp_unclustered} and \eqref{iter:u}.

The second term on the RHS of \eqref{eq:decomposed_dvdp_unclustered} sums over all the unclustered nodes $j\in\mathcal{N}_0$. It can be calculated using \eqref{derivesecond:p} by the CC and then sent to node $h$. In particular, $R_{ih}$ and $\mathds{1}(j\in\mathbb{P}_{h})$ are known to the CC for the calculation of \eqref{derivesecond:p}, since node $j$, its parent node $i$, and node $h$ are all in the backbone network.




To summarize, the computation of the key term $\left(\frac{\partial \boldsymbol{v}}{\partial \boldsymbol{u} }\right)^{\top} \left(\overline{\boldsymbol{\mu}}-\underline{\boldsymbol{\mu}}\right)$ in \eqref{iter:u} can be performed through the coordination of the CC and the RCs in a hierarchical manner. 
This inspires our design of Algorithm 1 to solve the OPF problem \eqref{OPF}.  
Compared to the conventional centralized primal-dual gradient method, the hierarchical Algorithm 1 can reduce computational complexity as analyzed in \cite{zhou2019accelerated}. Due to space limitation, the convergence proof of Algorithm 1 is deferred to the journal version of this work.    

\begin{algorithm*}
\caption{Hierarchical OPF Algorithm} \label{algCCRC}
\begin{algorithmic}[1]
\Repeat
\State At time step $t$, every node $h\in \mathcal{N}$ performs local update of primal variables:
\begin{subequations}
\begin{alignat}{2}
p_{h}(t+1)=\left[p_{h}(t)-\sigma_u\left(\frac{\partial f_{h}\left({p}_h(t), {q}_h(t)\right)}{\partial p_{h}}+\alpha_{h}(t)\right)  \right]_{\mathcal{Y}_{h}} \nonumber\\
q_{h}(t+1)=\left[q_{h}(t)-\sigma_u\left(\frac{\partial f_{h}\left({p}_h(t), {q}_h(t)\right)}{\partial q_{h}}+\beta_{h}(t)\right)  \right]_{\mathcal{Y}_{h}} \nonumber
\end{alignat}
\end{subequations}
where $\alpha_{h}(t)=\quad \sum_{j \in \mathcal{N}}\frac{\partial v_{j}\left(\boldsymbol{u}(t)\right)}{\partial p_{h}}\cdot\left(\overline{\mu}_{j}(t)-\underline{\mu}_{j}(t)\right)$ is given by \eqref{eq:decomposed_dvdp} if $h\in\mathcal{N}_k$ is in subtree $k$, or \eqref{eq:decomposed_dvdp_unclustered} if $h\in\mathcal{N}_0$ is unclustered; and $\beta_{h}(t)=\quad \sum_{j \in \mathcal{N}}\frac{\partial v_{j}\left(\boldsymbol{u}(t)\right)}{\partial q_{h}}\cdot\left(\overline{\mu}_{j}(t)-\underline{\mu}_{j}(t)\right)$. The updated decisions $\boldsymbol{u}(t+1)=\left(\boldsymbol{p}(t+1), \boldsymbol{q}(t+1)\right)$ are actuated by controllable devices. 

\State Based on local voltage measurement $v_h(t)$, every node $h\in\mathcal{N}$ updates its dual variables locally:
\begin{alignat}{2}
\underline{\mu}_{h}(t+1)=\left[\underline{\mu}_{h}(t)+\sigma_{\mu}\left(\underline{v}_{h}-v_{h}(t)-\epsilon \underline{\mu}_{h}(t)\right)\right]_{+}, \quad \overline{\mu}_{h}(t+1)=\left[\overline{\mu}_{h}(t)+\sigma_{\mu}\left(v_{h}(t)-\overline{v}_{h}-\epsilon \overline{\mu}_{h}(t)\right)\right]_{+}. \nonumber
\end{alignat}
\State Every RC $k\in\mathcal{K}$ calculates the following weighted sum of dual variables and sends it to the CC:
\begin{alignat}{2}
\sum_{j \in \mathcal{N}_{k}}\left(1-\frac{\left|z_{ij}\right|^{2}\ell_{ij}}{v_{i}}\right)\cdot\left(\overline{\mu}_{j}(t+1)-\underline{\mu}_{j}(t+1)\right). \nonumber
\end{alignat}
%
\State The CC computes the second and third terms on the RHS of \eqref{eq:decomposed_dvdp} for each destination node $h\in\mathcal{N}_k$ in a subtree $k$, and computes the first and second terms on the RHS of \eqref{eq:decomposed_dvdp_unclustered} for each unclustered destination node $h\in\mathcal{N}_0$. The CC adds up those terms and sends the result to the corresponding RC $k$ that hosts the destination node $h\in\mathcal{N}_k$, or to $h\in\mathcal{N}_0$ directly if it is unclustered.  
The terms needed for computing $\beta_h(t+1)$ are calculated and sent by the CC similarly.    

%
\State If $h \in \mathcal{N}_0$ is unclustered, it already obtains $\alpha_h(t+1)$ as \eqref{eq:decomposed_dvdp_unclustered}. Otherwise if $h\in\mathcal{N}_k$, RC $k$ computes the first term on the RHS of \eqref{eq:decomposed_dvdp}, adds it with the result received from the CC in step 5, and sends the result $\alpha_h(t+1)$ to the destination node $h$. The term $\beta_h(t+1)$ is obtained similarly. 



\Until $\left\|\boldsymbol{u}(t+1)-\boldsymbol{u}(t)\right\|_2 < \delta$ for some preset threshold $\delta>0$, or a maximum number of iterations is reached. 

\State Otherwise $t\leftarrow (t+1)$ and go back to step 2. 
\end{algorithmic}
\end{algorithm*}

\section{Numerical Experiments}\label{sec:simulation}
We conduct numerical experiments to demonstrate our improvement over the previous method \cite{zhou2019hierarchical} based on the linearized power flow model \eqref{simppowerequation}.

\subsection{Experiment setup}

\begin{figure} 
    \centering
    \includegraphics[width = 1.00\columnwidth]{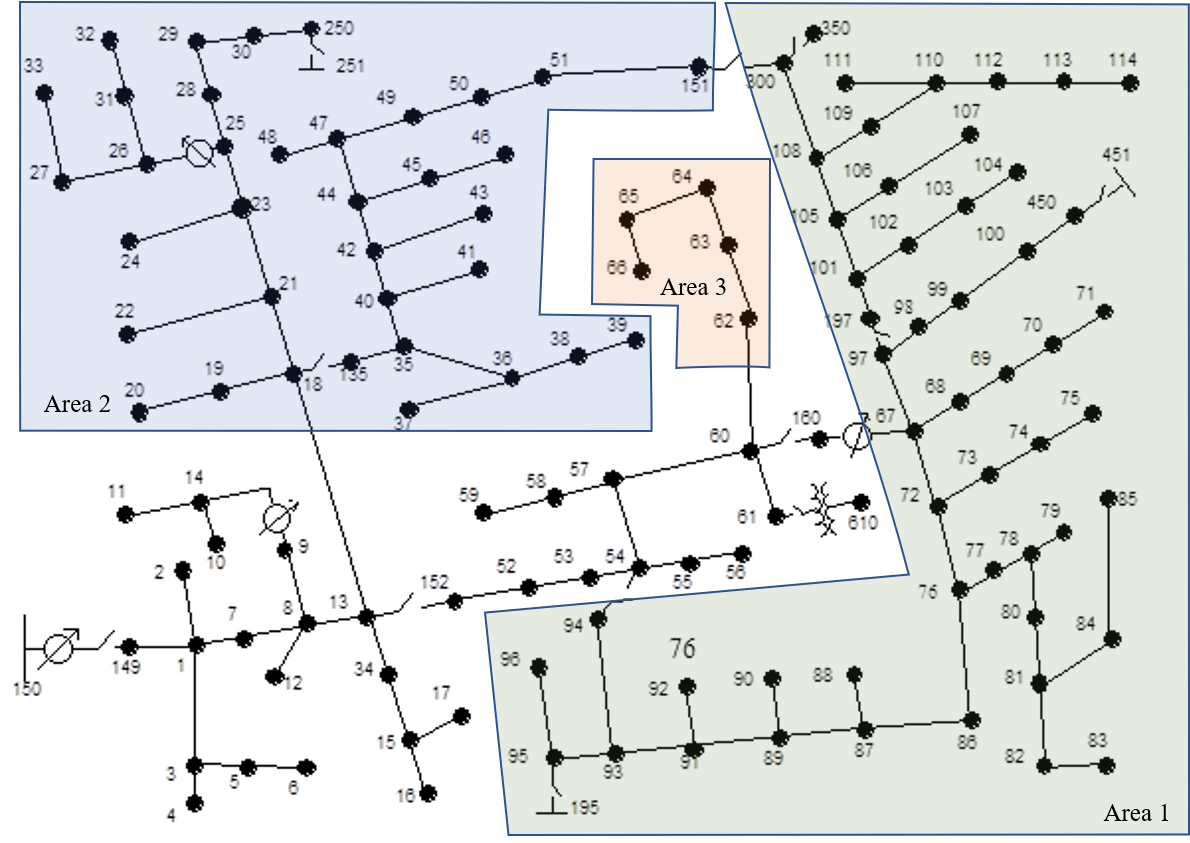}
    \caption{The clustering of IEEE 123-node network in our experiments.}
    \label{ieee123}
\end{figure}

We consider the single-phase equivalent models of the IEEE 37-node and 123-node test networks. They are clustered into subtrees as shown by Figure \ref{ieee37} (left panel) and Figure \ref{ieee123}, respectively. 

For the ease of illustration, we make the following modifications to the original network models on the IEEE PES website (https://cmte.ieee.org/pes-testfeeders/resources/): 
\begin{enumerate}
\item  We model each multi-phase line as a single-phase line with the average impedance of multiple phases. We also convert the multi-phase wye- and delta-connected loads into single-phase loads by taking their average over multiple phases. 
\item  We multiply the original load data of the 37-node network by six, and those of the 123-node network by two, to create scenarios with serious voltage issues.  
\item All the loads are treated as constant-power loads. Detailed models of capacitors, regulators, and breakers are not simulated. 
\end{enumerate}

For each network, let $\mathcal{N}_L$ denote the set of nodes that host nonzero loads.
We denote the \emph{negative net} power injections from the IEEE load data by $(\underline{p}_{h},\underline{q}_{h})$ for all load nodes $h\in\mathcal{N}_L$, and treat them as the nominal injections or the most preferred injections. 
The feasible power injection regions in \eqref{region} are defined as $\mathcal{Y}_h = \left\{(p_h,q_h) \mid \underline{p}_{h} \leq p_{h} \leq 0.3 \underline{p}_h<0, \underline{q}_{h} \leq q_{h} \leq 0.3\underline{q}_{h}<0\right\}$. 
The OPF objective function in \eqref{OPF:obj} is defined as $f_{h}(p_h,q_h)=(p_{h}-\overline{p}_{h})^{2}+(q_{h}-\overline{q}_{h})^{2}$ for all $h\in\mathcal{N}_L$, which aims to minimize the disutility caused by the deviation from the nominal loads. 

For each network, the voltage magnitude of the root (slack) node is fixed at $1.05$ per unit (p.u.). The lower and upper limits for safe voltage are set at $0.95$ p.u. and $1.05$ p.u., respectively. 
The step sizes for primal and dual variable updates are empirically chosen as $\sigma_u = 2\times 10^{-3}$ and $\sigma_{\mu}=1\times 10^{-3}$, respectively.

For each given power injection $\boldsymbol{u}=\left(\boldsymbol{p},\boldsymbol{q}\right)$, an OpenDSS power flow simulation is performed to obtain the corresponding voltage $\boldsymbol{v}(\boldsymbol{u})$ and the associated line-flow quantities $(\boldsymbol{P},\boldsymbol{Q}, \boldsymbol{\ell})$, which are treated as the feedback signals measured from the power network.  The Python 3.7 programs for OPF algorithms and the OpenDSS simulations are run on a laptop equipped with Intel Core i7-9750H CPU @ 2.6 GHz, 16 GB RAM, and Windows 10 Professional OS.

\subsection{Numerical results}


Due to space limitation, we only show the simulation results in the 123-node network model. The results in the 37-node network are similar.

\begin{figure}
\centering
\includegraphics[width=1.00\columnwidth]{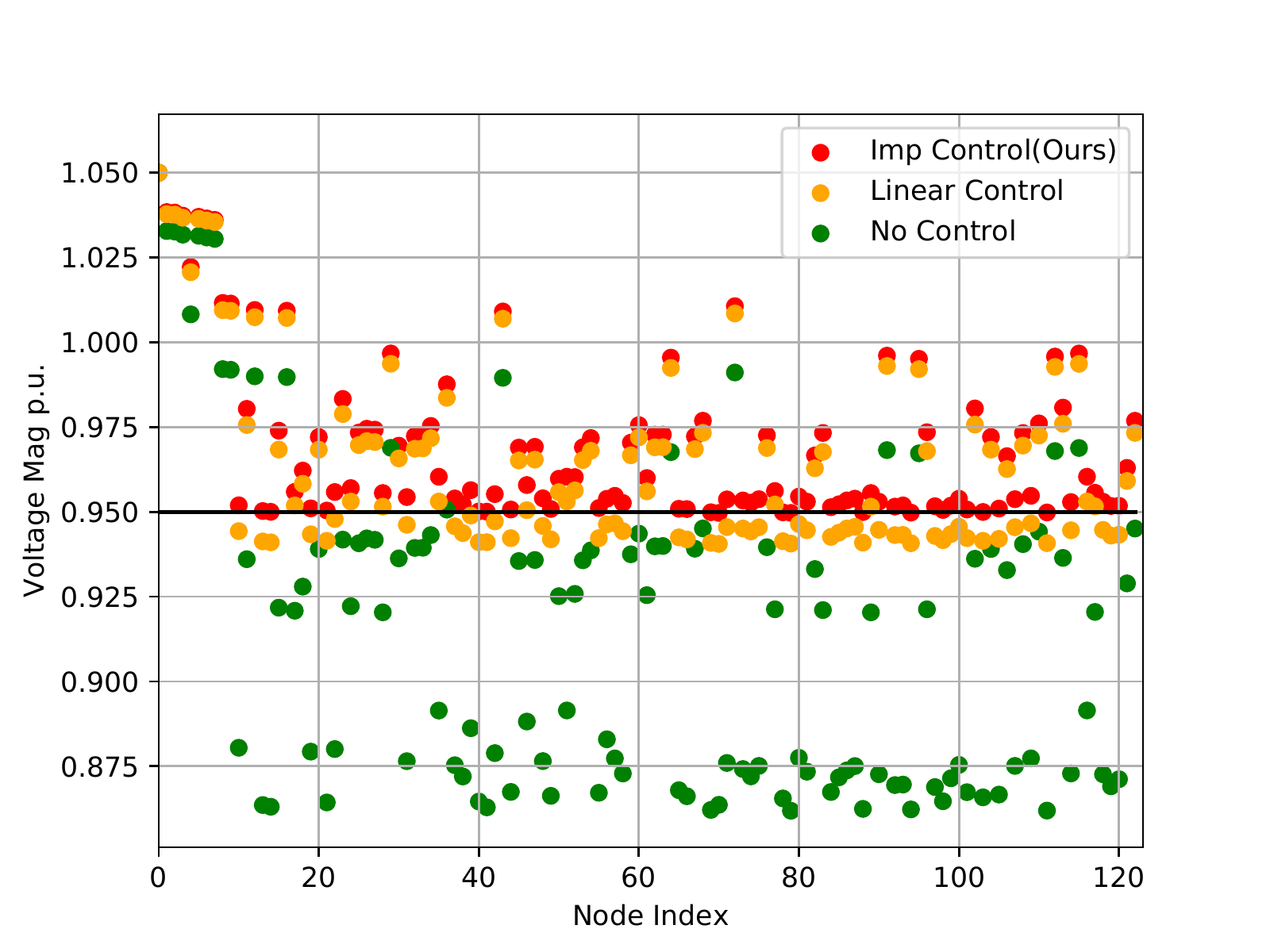}\label{fig:voltage_scatter_123}
\caption{The voltages in the IEEE 123-node network, in three cases: ``no control'' (using the nominal power injections); ``linear control'' (solving OPF with algorithm  \cite{zhou2019hierarchical} based on linear power flow model); ``improved control'' (based on the proposed Algorithm 1 with improved gradient evaluation).}\label{fig:voltage_scatter}
\end{figure} 

\textit{Voltage safety}: The voltage magnitudes at different nodes of the 123-node network are plotted in Figure \ref{fig:voltage_scatter}, for three cases. The first case, referred to as ``no control'', just takes the nominal power injections $(\underline {\boldsymbol{p}},\underline{\boldsymbol{q}})$ from the IEEE data, which causes severe violation of the voltage lower bound. The second case, referred to as ``linear control'', applies the primal-dual gradient algorithm in \cite{zhou2019hierarchical} based on the linearized power flow model \eqref{simppowerequation} to determine power injections, which enhances the overall voltage level but still leaves most nodes below the voltage lower bound. The third case, referred to as ``improved control (ours)'', solves the OPF with the proposed Algorithm 1 based on our improved gradient evaluation method, which effectively lifts all the voltages above the lower bound. This result verifies our analysis in Section \ref{sec:method} that the improved gradient evaluation can prevent voltage violation caused by the linear model that optimistically estimates the voltages to be safe.

\begin{figure}
\centering
\includegraphics[width=1.00\columnwidth]{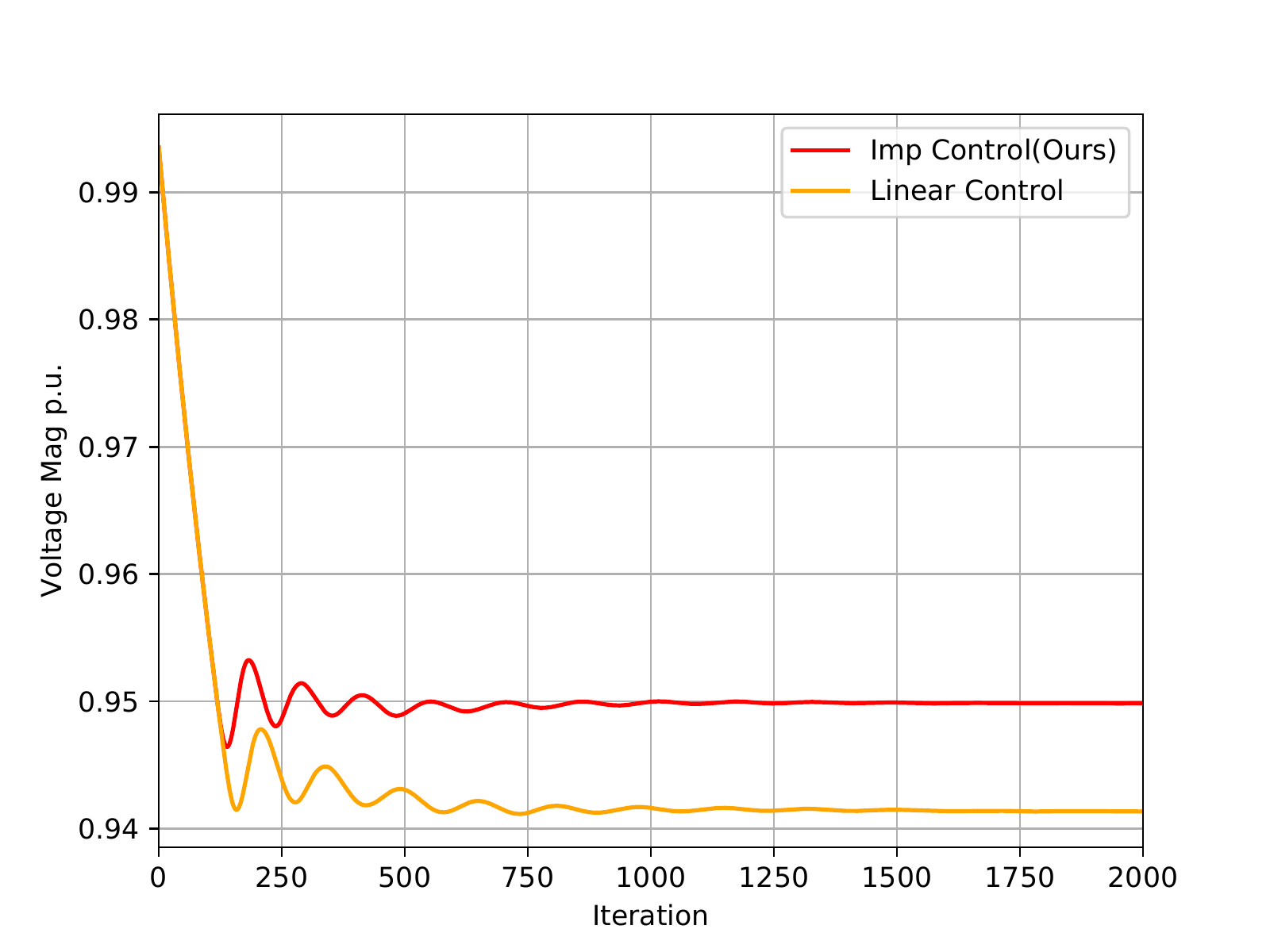}
\caption{The change of voltage at a node of the IEEE 123-node network, in the proposed Algorithm 1 (``improved control'') and the OPF algorithm in \cite{zhou2019hierarchical} based on the linear power flow model (``linear control'').}\label{fig:voltage_convergence_123}
\end{figure} 

\textit{Computational efficiency}: Figure \ref{fig:voltage_convergence_123} shows the change of voltage at a node of the 123-node network (similar trends are displayed for the voltages at other nodes). This result verifies again that the voltage converges to a safe level under the proposed Algorithm 1, compared to the previous algorithm \cite{zhou2019hierarchical} that leads to a violation of voltage lower bound. 
Moreover, Algorithm 1 converges in fewer iterations. 
Indeed, Algorithm 1 spends 205 seconds to complete 2,000 iterations, compared to the previous algorithm in \cite{zhou2019hierarchical} which spends 166 seconds. This would be an acceptable increase, considering the improved voltage safety achieved by Algorithm 1. 

\section{Conclusion}\label{sec:conclusion}

We proposed a hierarchical OPF algorithm based on a new gradient evaluation method that is more accurate than previous methods based on the linearized power flow model. The proposed method can prevent the previous voltage violation with moderate extra computations, as verified by numerical results on IEEE networks. 

In future work, we will provide formal convergence proof and complexity analysis of the proposed algorithm. We will also extend the proposed algorithm to multi-phase networks, as what was done from \cite{zhou2019hierarchical} to
 \cite{zhou2019accelerated}.

 {\appendix[Proof of Lemma 1]
\begin{figure}
\centering
\includegraphics[width=0.90\columnwidth]{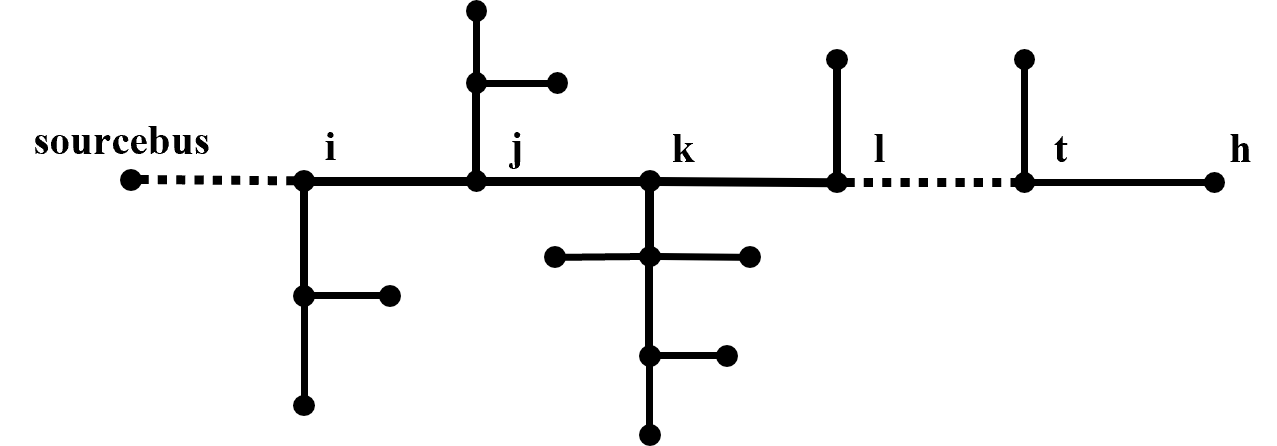}
\caption{Index the nodes in a tree graph.} \label{fig:treegraph}
\end{figure} 
\textit{Proof:} Take the difference between the linear model (\ref{simppowerequation}) and the actual nonlinear model (\ref{powerequation}):
\begin{subequations}
\begin{alignat}{2} 
\hat{v}_{j}-v_{j}=&\hat{v}_{i}-v_{i}-2\left(r_{ij}\left(\hat{P}_{ij}-P_{i j}\right)+x_{ij}\left(\hat{Q}_{ij}- Q_{i j}\right)\right)\nonumber\\
&-(r_{ij}^{2}+x_{ij}^{2})\ell_{ij},  \label{error::v} \\
\hat{P}_{ij}-P_{i j}&=\sum_{k:(j,k)\in \mathcal{E}}(\hat{P}_{jk}-P_{j k})-r_{ij}\ell_{ij},  \label{error::P}\\
\hat{Q}_{ij}-Q_{i j}&=\sum_{k:(j,k)\in \mathcal{E}}(\hat{Q}_{jk}-Q_{j k})-x_{ij}\ell_{ij}.  \label{error::Q}
\end{alignat}
\end{subequations}
Substitute (\ref{error::P})--(\ref{error::Q}) into (\ref{error::v}), which yields:
\begin{alignat}{2} \label{error::v:2}
\hat{v}_j-&v_j =\hat{v}_i-v_i-2\left(r_{i j} \sum_{k:(j, k) \in \mathcal{E}}\left(\hat{P}_{j k}-P_{j k}\right)\right. \nonumber\\
&\left.+x_{i j} \sum_{k:(j, k) \in \mathcal{E}}\left(\hat{Q}_{j k}-Q_{j k}\right)\right)+\left(r_{i j}^2+x_{i j}^2\right) \ell_{i j}.
\end{alignat}

We further calculate the active and reactive power loss terms $\sum_{k:(j, k) \in \mathcal{E}}(\hat{P}_{jk}-P_{jk})$ and $\sum_{k:(j, k) \in \mathcal{E}}(\hat{Q}_{j k}- Q_{j k})$ on the RHS of (\ref{error::v:2}). Following the tree structure of the network (i.e., a tree graph shown by Figure \ref{fig:treegraph}), we recursively calculate (\ref{error::P})--(\ref{error::Q}) until reaching leaves:
\begin{subequations}
\begin{alignat}{2}
\hat{P}_{jk}-P_{jk}&=\sum_{l:(k,l)\in \mathcal{E}}(\hat{P}_{kl}-P_{kl})-r_{jk}\ell_{jk},  \label{error::P:jk}\\
\hat{Q}_{jk}-Q_{jk}&=\sum_{l:(k,l)\in \mathcal{E}}(\hat{Q}_{kl}-Q_{kl})-x_{jk}\ell_{jk},  \label{error::Q:jk} \\
\qquad \dots\nonumber\\
\hat{P}_{th}-P_{t h}&=-r_{t h}\ell_{t h}  \label{error::P:th}, ~\text{for leaf nodes $h$},\\
\hat{Q}_{th}-Q_{t h}&=-x_{t h}\ell_{t h}  \label{error::Q:th}, ~\text{for leaf nodes $h$}.
\end{alignat}
\end{subequations}

We substitute the above equations layer-by-layer, i.e., (\ref{error::P:th}), $\dots$, into (\ref{error::P:jk}) for active power and (\ref{error::Q:th}), $\dots$, into (\ref{error::Q:jk}) for reactive power. This gives the expression of the active and reactive power loss terms as:
\begin{subequations} \label{error::jk:3}
\begin{alignat}{2}
&\sum_{k:(j, k) \in \mathcal{E}}(\hat{P}_{j k}-P_{j k})= \\
&\qquad -\sum_{k:(j,k)\in \mathcal{E}}\sum_{l:(k,l)\in \mathcal{E}}\dots\sum_{h:(t,h)\in \mathcal{E}}r_{th}\ell_{th}-\dots - \nonumber \\
&\qquad -\sum_{k:(j,k)\in \mathcal{E}}\sum_{l:(k,l)\in \mathcal{E}}r_{kl}\ell_{kl}-\sum_{k:(j,k)\in \mathcal{E}}r_{jk}\ell_{jk},  \label{error::jk:P:3}\\
&\sum_{k:(j, k) \in \mathcal{E}}(\hat{Q}_{j k}-Q_{j k})= \\
&\qquad -\sum_{k:(j,k)\in \mathcal{E}}\sum_{l:(k,l)\in \mathcal{E}}\dots\sum_{h:(t,h)\in \mathcal{E}}x_{th}\ell_{th} -\dots - \nonumber \\
&\qquad -\sum_{k:(j,k)\in \mathcal{E}}\sum_{l:(k,l)\in \mathcal{E}}x_{kl}\ell_{kl}-\sum_{k:(j,k)\in \mathcal{E}}x_{jk}\ell_{jk}. \label{error::jk:Q:3}
\end{alignat}
\end{subequations}

We define $g_{j}$ and $\eta_{j}$ as the \textit{additive inverse} of RHS of (\ref{error::jk:P:3}) and (\ref{error::jk:Q:3}), respectively. In fact, $g_{j}$ and $\eta_{j}$ sum the active and reactive power loss of the lines downstream of node $j$. We have more compact expressions for $g_{j}$ and $\eta_{j}$:
\begin{subequations} \label{g:eta:2}
\begin{alignat}{2} \nonumber
g_{j}=\sum_{(\alpha,\beta)\in \operatorname{down}(j)}r_{\alpha\beta}\ell_{\alpha\beta}, \quad \eta_{j}=\sum_{(\alpha,\beta)\in \operatorname{down}(j)}x_{\alpha\beta}\ell_{\alpha\beta},
\end{alignat}
\end{subequations}
where $\operatorname{down}(j)$ denotes the set of lines downstream of node $j$. From the tree structure, we have:
\begin{subequations} \label{g:eta:3}
\begin{alignat}{2}
g_{j}=\sum_{k:(j,k)\in \mathcal{E}}r_{jk}\ell_{jk}+\sum_{k:(j,k)\in \mathcal{E}}g_{k},\\
\eta_{j}=\sum_{k:(j,k)\in \mathcal{E}}x_{jk}\ell_{jk}+\sum_{k:(j,k) \in \mathcal{E}}\eta_{k}.
\end{alignat}
\end{subequations}

The expressions of (\ref{g:eta:3}) imply that a backward sweep from the leaves to the root node in a power distribution network will give all the results of $(g_{j},\eta_{j}), \forall j \in \mathcal{N}$. With $(g_{j},\eta_{j})$, the voltage error (\ref{error::v:2}) becomes:
\begin{alignat}{2} \label{error::v4}
\hat{v}_{j}-v_{j}&=\hat{v}_{i}-v_{i}+2\left(r_{ij}g_{j}+x_{ij}\eta_{j}\right)+(r_{ij}^{2}+x_{ij}^{2})\ell_{ij}.
\end{alignat}

Equation (\ref{error::v4}) implies that the voltage error is accumulated from the root to the leaves. A forward sweep from the root to leaves will return voltage error for all nodes $h\in \mathcal{N}$ as:

\begin{alignat}{2} \label{verror::2}
\hat{v}_{h}-v_{h}= \sum_{(\zeta,\xi)\in \mathbb{P}_{h}}\left(2\left(r_{\zeta\xi}g_{\xi}+x_{\zeta\xi}\eta_{\xi}\right)+(r_{\zeta\xi}^2+x_{\zeta\xi}^2)\ell_{\zeta\xi}\right). \nonumber
\end{alignat}

This completes the proof of \emph{Lemma \ref{lemma::1}.}

\bibliographystyle{ieeetr}
\bibliography{ref}

\end{document}